\documentclass{amsart}

\begin{document}

\title{Polar orbits around binary stars}
\author{Greg Egan}
\email{gregegan@netspace.net.au}

\date{Dec 29, 2017}

\begin{abstract} 
\noindent 
Oks proposes the existence of a new class of stable planetary orbits around binary stars, in the shape of a helix on a conical surface whose axis of symmetry coincides with the interstellar axis.  We show that this claim relies on the inappropriate use of an effective potential that is only applicable when the stars are held motionless, and that the existence of the torques required to maintain the proposed orbits would have been empirically detected in the motion of artificial satellites in high-inclination orbits.
\end{abstract}

\maketitle

\section{Introduction}
The restricted three-body problem, in which two bodies execute a Kepler orbit while a third body of negligible mass moves in the resulting potential, has been the subject of intensive investigation since the time of Newton, with important contributions by Euler, Jacobi, Lagrange and Poincar\'{e} \cite{Murray}.  Given this long history, it would be remarkable if an entirely new class of stable solutions were to be discovered in the twenty-first century.

Oks \cite{Oks2015} proposes the existence of a new class of solutions, in which the body of negligible mass orbits, approximately, in a plane that remains orthogonal to the axis between the two massive bodies as they move in their own orbit, while undergoing perturbations on the surface of a cone that is aligned with the same axis. 

The existence of such orbits would be extraordinary, not only by virtue of being a novel result in a well-studied field, but because they would exhibit a strongly counter-intuitive effect:  the  gyroscopic tendency of an orbiting body to resist the rotation of its orbital plane would need to be overcome by a torque from the more massive bodies, and it is not at all obvious how such a torque would arise.

Our analysis, however, does not support the claims in \cite{Oks2015}.  Furthermore, we find that the claims imply empirically false predictions about artificial satellites in high-inclination orbits.

\section{Analysis in the original paper}

Oks \cite{Oks2015} begins by considering planetary orbits around the interstellar axis between two \emph{motionless} stars.  Such systems can be characterized by a two-dimensional parameter space $(w,b)$, where $w$ is the scaled distance of the plane of the orbit along the interstellar axis from one of the stars, and $b$ is the mass ratio of the stars.  An effective potential $U_{\mathrm{eff}}$ is derived, associated with a conserved quantity $M$, the projection of the planet's orbital angular momentum on the interstellar axis, and the regions of the parameter space that allow points of stable equilibrium for $U_{\mathrm{eff}}$ are identified.

Further restrictions on $(w,b)$ are derived which force the frequency of the planetary orbit to be much faster than the Kepler frequency of the binary pair. The aim of these restrictions is to allow a separation of the system into rapid and slow subsystems, but while the paper identifies a range of orbits that meet this condition on the orbital frequencies, it neither cites nor proves any result to the effect that such a separation is a valid method of approximation in this particular context.

The problem of the orbit when the binary pair rotates is studied in a non-inertial frame in which the stars remain fixed and the dynamics is modified by centrifugal and Coriolis forces.  This analysis proceeds by treating the planet's degrees of freedom in a plane through the interstellar axis as a two-dimensional harmonic oscillator, driven by the Coriolis force, around a stable equilibrium point in the $U_{\mathrm{eff}}$ derived under the assumption that the stars are motionless.  If the amplitudes of the driven oscillator are small, the full motion of the planet will be a helix lying on the surface of a cone whose axis of symmetry is the interstellar axis.

However, we believe this analysis is not correct, as $U_{\mathrm{eff}}$ depends on $M$, which is not conserved, even approximately, as the binary pair rotates.

\section{Failure of $M$ to be conserved}
In Section 3.2 of \cite{Oks2015}, the centrifugal and Coriolis forces are obtained by approximating the planetary orbit as circular, with its centre on the interstellar axis and its plane orthogonal to that axis.  This is not intended as a detailed description of the motion, but the departures from it are assumed to be small and cyclical, with periods much shorter than the time scale associated with the stellar orbit.

Our analysis will extend that in \cite{Oks2015} by computing the torques acting on the planet, and rather than assuming that the plane of the orbit remains orthogonal to the axis, we will allow for the possibility that its normal vector makes an angle $\alpha(t)$ with the axis, while remaining in the plane of the stellar orbit.  We will also allow for the possibility that the centre of the orbit moves along the axis.

We will work in Cartesian coordinates in a rotating frame fixed to the stars, with the origin at the centre of mass of the binary pair, the $x$-axis as the interstellar axis, and the $z$-axis as the axis of rotation for the stellar orbit.  Wlog, we will assume that the planet has unit mass.

$M(t)$, the projection of the planet's angular momentum on the interstellar axis, is given by:

\begin{equation}
\label{E:M}
M(t) = y(t) z'(t) - z(t) y'(t)
\end{equation}

If the angular velocity of the stellar orbit is $\omega$, the Coriolis force is:

\begin{equation}
\label{E:F1}
F_1(t) = 2 \omega \:  (y'(t), -x'(t), 0)
\end{equation}
 
The component of the torque along the $x$-axis due to the Coriolis force is:

\begin{equation}
\label{E:T1}
\begin{array}{rcl}
T_1(t) & = & y(t) F_{1,z}(t) - z(t) F_{1,y}(t) \\
& = & 2 \omega  z(t) x'(t)
\end{array}
\end{equation}

The centrifugal force is:

\begin{equation}
\label{E:F2}
F_2(t) = \omega^2 \: (x(t), y(t), 0)
\end{equation}

The component of the torque along the $x$-axis due to the centrifugal force is:

\begin{equation}
\label{E:T2}
\begin{array}{rcl}
T_2(t) & = & y(t) F_{2,z}(t) - z(t) F_{2,y}(t) \\
& = & -\omega^2 y(t) z(t)
\end{array}
\end{equation}

The gravitational force exerted by either star, the displacement of the planet from the centre of mass, and the interstellar axis all lie in the same plane, so the torque due to the gravitational force is always orthogonal to the interstellar axis, and has no effect on $M(t)$.

We approximate the motion of the planet with a circular orbit of angular velocity $f$ and radius $\rho$, centred on the point $(x_1(t), 0, 0)$, and tilted at an angle $\alpha(t)$.

\begin{equation}
\label{E:xyz}
\begin{array}{rcl}
x(t) & = & x_1(t) + \rho  \sin(\alpha(t) ) \cos(f t) \\
y(t) & = & \rho  \cos(\alpha(t) ) \cos(f t) \\
z(t) & = & \rho  \sin(f t)
\end{array}
\end{equation}

Substituting (\ref{E:xyz}) into (\ref{E:M}), (\ref{E:T1}) and (\ref{E:T2}) we obtain:

\begin{equation}
\label{E:MC}
M(t)  = \frac{1}{2} \rho ^2 \sin (2 f t) \sin (\alpha (t)) \alpha '(t) + f \rho ^2 \cos (\alpha (t))
\end{equation}

\begin{equation}
\label{E:T1C}
\begin{array}{rcl}
T_1(t) &= & \rho ^2 \omega  \sin (2 f t) \cos (\alpha (t)) \alpha '(t) +f \rho ^2 \omega  \cos (2 f t) \sin (\alpha (t)) \\
& & -\rho  \omega  \left(f \rho  \sin (\alpha (t))-2 \sin (f t) x_1'(t)\right)
\end{array}
\end{equation}

\begin{equation}
\label{E:T2C}
T_2(t) =-\frac{1}{2} \rho ^2 \omega ^2 \sin (2 f t) \cos (\alpha (t))
\end{equation}

Differentiating (\ref{E:MC}), we obtain:

\begin{equation}
\label{E:MCD}
\begin{array}{rcl}
M'(t)  &=& f \rho ^2 \sin (\alpha (t)) \alpha '(t) (\cos (2 f t) - 1) \\
& & +\frac{1}{2} \rho ^2 \sin (2 f t) \left(\sin (\alpha (t)) \alpha ''(t)  + \cos (\alpha (t)) \alpha '(t)^2 \right)
\end{array}
\end{equation}

The rate of change of the planet's angular momentum must be equal to the total torque it experiences, so for the components projected on the interstellar axis we have:

\begin{equation}
\label{E:MTT}
M'(t) - T_1(t)  - T_2(t) = 0 
\end{equation}

If we substitute (\ref{E:T1C}),  (\ref{E:T2C}) and  (\ref{E:MCD}) into (\ref{E:MTT}), and require the terms with factors of $\sin(2 f t)$, $\cos(2 f t)$ and $\sin(f t)$ to be equal to zero, we obtain the differential equations:

\begin{equation}
\label{E:DIFF}
\begin{array}{rcl}
\alpha ''(t) \sin (\alpha (t))+\cos (\alpha (t)) \left(\omega -\alpha '(t)\right)^2 & = & 0 \\
\sin (\alpha (t)) \left(\omega -\alpha '(t)\right) & = & 0 \\
x_1'(t) & = & 0 \\
\end{array}
\end{equation}

The second of these equations is satisfied by either $\alpha(t)=n \pi$ for some integer $n$, or $\alpha(t) = \omega t + \alpha_0$, but the first equation is satisfied only by the latter choice.  So the only solution to the full system of equations is:

\begin{equation}
\label{E:SOL}
\begin{array}{rcl}
\alpha(t) & = & \omega t + \alpha_0 \\
x_1(t) & = & x_0
\end{array}
\end{equation}

If we substitute this solution into (\ref{E:T1C}),  (\ref{E:T2C}) and  (\ref{E:MCD}), then equation (\ref{E:MTT}) is satisfied, because the only remaining terms that we did not use to construct (\ref{E:DIFF}) are a multiple of the left-hand side of the second equation in (\ref{E:DIFF}).

Substituting this solution into (\ref{E:MC}), we obtain:

\begin{equation}
\label{E:MCSOL}
M(t)  = f \rho ^2 \cos \left(\omega t + \alpha_0 \right) + \frac{1}{2} \rho ^2 \omega  \sin (2 f t) \sin \left(\omega t + \alpha_0 \right)
\end{equation}

Clearly $M(t)$ is not conserved on the time scale of the stellar orbit.  The regime of interest in \cite{Oks2015} is $f\gg \omega$, so the first term here will dominate, and $M(t)$ will be approximately sinusoidal, with a period equal to that of the stellar orbit.

This result shows that the treatment of the orbits in Section 3.2 of  \cite{Oks2015}, which relies on $M(t)$ being at least approximately conserved in order that the same effective potential can be used throughout the stellar orbit, is not correct.

\section{Empirical Implications}

The claims in \cite{Oks2015} can be tested by comparing their predictions to the behaviour of known systems.  While no extrasolar planets have been observed in suitably configured orbits, we can take the Earth and the Sun to be the two massive bodies, and a satellite orbiting the Earth with a high inclination as the third body.  While such a system is further subject to the effects of the Earth's oblateness, we can still ask whether the torques implied by \cite{Oks2015} have been observed at all, even in the presence of additional effects.

One complication is that both \cite{Oks2015} and a subsequent erratum \cite{Oks2016} contain errors in the algebra leading up to the final statements of the claimed perturbations from the equilibrium orbits, as given by equation (51) in both versions.  The original paper omitted a factor of twice the radius of the planetary orbit, and while the erratum correctly inserts this factor, it also corrects a sign error in the change of coordinates from equation (48) to equation (50), but does \emph{not} correct a sign error in the second line of equation (48), which was previously cancelling out the first error.  The result is that the erratum gives an incorrect version of equation (51) that claims the perturbations from the equilibrium would be entirely axial.  This makes no sense physically, since it amounts to saying that the response of a two-dimensional driven harmonic oscillator is colinear with the driving force (the Coriolis force, which is directed along the interstellar axis) even though the interstellar axis is \emph{not} a principal axis of the oscillator.

A correct calculation of the amplitudes of the perturbations yields a result that agrees with equation (51) in the original paper, except for the insertion of the missing factor of twice the radius of the planetary orbit:

\begin{equation}
\label{E:PERT}
\begin{array}{rcl}
\delta w(\tau) & = & \frac{4  \omega_s f_p v_0(w,b)}{\omega_+^2-\omega_-^2}  \cos (2 \alpha ) \cos (f_p \tau)\\
\delta v(\tau) & = & \frac{4 \omega_s f_p v_0(w,b)}{\omega_+^2-\omega_-^2}  \sin (2 \alpha ) \cos (f_p \tau)
\end{array}
\end{equation}

Here we have used the notation of \cite{Oks2015}, where $w$ and $v$ are scaled axial and radial coordinates in a cylindrical coordinate system whose axis is the interstellar axis, $v_0(w,b)$ is the scaled radius of the equilibrium orbit, $\tau$ is a scaled time coordinate, $f_p$ is the scaled primary frequency of the planetary motion, $\omega_s$ is the scaled Kepler frequency of the stars' orbit, $\omega_+$ and $\omega_-$ are scaled eigenfrequencies of the two-dimensional harmonic oscillator that approximates the well in the effective potential, and $\alpha$ is the angle between the coordinate directions and the principal axes of the harmonic oscillator. The scaling corresponds to a choice of units where the distance between the stars, the gravitational constant $G$, and the mass of the lighter of the two stars, are all set equal to 1.

The Earth-Sun system has a mass ratio, $b$, of approximately $3.33 \times 10^5$.  We will choose $w = 5 \times 10^{-8}$, which describes an orbit that lies roughly in a plane orthogonal to the Earth-Sun axis and displaced by $5 \times 10^{-8}$ au, or about 7.5 km, from the centre of the Earth towards the sun.

This choice meets all the criteria for the analysis in \cite{Oks2015} to apply.  Equation (12) of \cite{Oks2015} yields a scaled equilibrium orbital radius of $v_0(w,b)=5.315 \times 10^{-5}$, which corresponds to about 7972 km.  The equilibrium in the effective potential here is stable, since $v_0(w,b)$ exceeds the critical value $v_{\mathrm{crit}}(w)$, defined in inequality (26) of \cite{Oks2015}, which is $1.412 \times 10^{-7}$.  The frequency of an orbit around the Earth at this radius is about 4470 times greater than the frequency of the Earth's orbit around the Sun.  And finally, the amplitudes of the scaled axial and radial perturbations, given by (\ref{E:PERT}) in the present paper, are $8.37 \times 10^{-9}$ and $8.42 \times 10^{-6}$, which correspond to 1.25 and 1263 km respectively.  While the radial perturbations are not negligible compared to the total radius of the orbit, they lie within the region of the effective potential where it is well-approximated as a two-dimensional harmonic oscillator.

If the analysis in  \cite{Oks2015} were correct, any satellite's orbit with a plane roughly orthogonal to the Earth-Sun axis and a mean radius of about 8000 km would experience a  rotation relative to the fixed stars, keeping it facing towards the Sun.  While such orbits, known as dawn/dusk Sun-synchronous orbits, are desirable for solar observatories, in fact they are only achieved by exploiting the oblateness of the Earth \cite{Capderou}.  If a torque were available to provide the same effect independently of the shape of the Earth, it would be of the same order of magnitude as the torque due to oblateness, so its existence would have been observed and documented by now.

\section{Conclusions}

The claims in \cite{Oks2015} in regard to a proposed new class of stable orbits in the restricted three-body problem are based on an incorrect assumption:  that the effective potential that governs a planet orbiting the axis between two motionless stars can still be employed when the stars are rotating.  An analysis of the torques acting on the planet shows that its orbital angular momentum will change in a manner that is inconsistent with this assumption. 

Furthermore, if the claims were correct, their empirical consequences would have been manifest in the detailed behaviour of high-inclination satellite orbits, and this is not the case.

\section*{Acknowledgement}

This is a preprint.  The final version of the paper is: Egan, G. Celest Mech Dyn Astr (2018) 130: 5. https://doi.org/10.1007/s10569-017-9803-7

\end{document}